\documentclass[]{article}
\begin{document}
\newtheorem{proposition}{Proposition}[section]
\newtheorem{definition}{Definition}[section]
\newtheorem{lemma}{Lemma}[section]

\title{\bf A Class of Noncommutative Spectra}
\author{Keqin Liu\\Department of Mathematics\\The University of British Columbia\\Vancouver, BC\\
Canada, V6T 1Z2}
\date{March, 2011}
\maketitle

\begin{abstract} We construct a class of noncommutative spectra and give the basic properties of the class of noncommutative spectra.
\end{abstract}

How to select a class of noncommutative rings such that we can rewrite algebraic geometry in the context of the class of  noncommutative rings? The purpose of this paper is to present our answer to this interesting problem. The main result of this paper is to construct  a class of noncommutative spectra which are called trispectra. Trispectra come from a class of nonncommutative rings which are called Hu-Liu trirings by us. In section 1, we introduce Hu-Liu trirings and give the basic properties of Hu-Liu trirings. The most important example of Hu-Liu trirings is triquaternions which is defined in section 1. Triquaternions, which are regarded as a kind of new numbers by us,  can be used to replace complex numbers to develop the counterpart of complex algebraic geometry. In section 2, we introduce prime triideals and prove some basic facts about prime triideals. In section 3, we use prime triideals to characterize the trinilradical. In the last section of this paper, we make trispectra into a topological space by introducing extended Zariski topology. 

\medskip
Throughout this paper, the word ``ring'' means an associative ring with an identity. A ring $R$ is also denoted by $(R,\, +,\, \cdot)$ to indicate that $+$ is the addition and $\cdot$ is the multiplication in the ring $R$. The word ``triring'' always means ``Hu-Liu triring''.

\bigskip
\section{Basic Definitions}

Let $A$ and $B$ be two subsets of a ring $(\, R, \, + , \, \cdot\,)$. We shall use $A+B$ and $AB$ to denote the following subsets of $R$ 
$$A+B:=\{\, a+b \, | \, a\in A, \, b\in B \, \}, \quad 
AB:=\{\, ab \, |\, a\in A, \, b\in B \,\}.$$

\medskip
Hu-Liu trirings are a class of noncommutative rings which can be used to rewrite the theory of commutative rings and algebraic geometry over commutative rings.

\medskip
\begin{definition}\label{def2.1} A ring $R$ with a multiplication $\cdot$ is called a {\bf Hu-Liu triring} if the following three properties hold.
\begin{description}
\item[(i)] There exist two commutative subrings $R_0$ and $R_1$ of the ring $(R_0, +, \cdot)$, called the {\bf even part} and {\bf odd part} of $R$ respectively, such that 
$R=R_0\oplus R_1$ (as Abelian groups) and 
\begin{equation}\label{eq1.1}
R_0R_0\subseteq R_0,\quad R_0R_1+R_1R_0\subseteq R_1, \quad R_1R_1=0;
\end{equation}
\item[(ii)] There exists a  binary operation $\sharp$ on the odd part $R_1$ such that 
$(\,R_1 , \, + , \, \,\sharp\, \,)$ is a commutative ring and the two associative products $\cdot$ and $\,\sharp\,$ satisfy the {\bf triassociative law}:
\begin{eqnarray}
\label{eq2} x (\alpha \,\sharp\, \beta) &=&(x \alpha) \,\sharp\, \beta, \\
\label{eq3} (\alpha \,\sharp\, \beta ) x &=&\alpha \,\sharp\, (\beta x), 
\end{eqnarray}
where $x\in R$ and $\alpha$, $\beta\in R_1$;
\item[(iii)] For each $x_0\in R_0$, we have
\begin{equation}\label{eq6}
R_1x_0=x_0R_1.
\end{equation}
\end{description}
\end{definition}

\medskip
A Hu-Liu triring $R=R_0\oplus R_1$ is sometimes denoted by \linebreak
$ (\,R=R_0\oplus R_1 , \, + , \, \cdot ,\,  \,\sharp\, \,)$,
where the associative product $\sharp$ on the odd part $R_1$ is called the {\bf local product}, and the identity $1^\sharp$ of the ring $(\, R_1 , \, + , \, \,\sharp\,\,)$ is called the {\bf local identity} of the triring $R$.

\medskip
Clearly, a commutative ring is a Hu-Liu triring with zero odd part. The first and the most important example of Hu-Liu trirings which is not a commutative ring is triquaternions whose definition is given in the following example.

\medskip
\noindent
{\bf Example} Let 
${\bf Q}=\mathcal{R}\,1\oplus \mathcal{R}\,i\oplus\mathcal{R}\,j\oplus\mathcal{R}\,k$ be a 
$4$-dimensional real vector space, where $\mathcal{R}$ is the field of real numbers. Then 
$(\,{\bf Q}={\bf Q}_0\oplus {\bf Q}_1, \, +, \,\cdot, \, \,\sharp\, \,)$ is a Hu-liu triring, where ${\bf Q}_0=\mathcal{R}\,1\oplus \mathcal{R}\,i$ is the even part, 
${\bf Q}_1=\mathcal{R}\,j\oplus\mathcal{R}\,k$ is the odd part, the ring multiplication $\cdot$ and the local product $\sharp$ are defined by the following multiplication tables:
\begin{center}
{\huge\begin{tabular}{|c||c|c|c|c|}
\hline
$\cdot$&$1$&$i$&$j$&$k$\\
\hline\hline
$1$&$1$&$i$&$j$&$k$\\
\hline
$i$&$i$&$-1$&$k$&$-j$\\
\hline
$j$&$j$&$-k$&$0$&$0$\\
\hline
$k$&$k$&$j$&$0$&$0$\\
\hline
\end{tabular}}
\qquad\qquad
{\huge\begin{tabular}{|c||c|c|}
\hline
$\sharp$&$j$&$k$\\
\hline\hline
$j$&$j$&$k$\\
\hline
$k$&$k$&$-j$\\
\hline
\end{tabular}}
\end{center}
The Hu-Liu triring ${\bf Q}$ is called the {\bf triquaternions}.

\hfill\raisebox{1mm}{\framebox[2mm]{}}

\medskip
{\bf Remark} For convenience, a Hu-Liu triring will be simply called a triring in the remaining part of this paper.

\medskip
\begin{definition}\label{def2.2} Let $(\,R=R_0\oplus R_1, \, +, \,\cdot, \, \,\sharp\, \,)$ be a triring, and let $I$ be a subgroup of the additive group of $R$.
\begin{description}
\item[(i)] $I$ is called a {\bf triideal} of $R$ if $IR+RI\subseteq I$, 
$I=(R_0\cap I)\oplus (I\cap R_1)$, and $I\cap R_1$ is an ideal of the ring 
$(\,R_1, \, +, \, \,\sharp\, \,)$.
\item[(ii)] $I$ is called a {\bf subtriring} of $R$ if $1\in I$, $II\subseteq I$,  
$I=(R_0\cap I)\oplus (I\cap R_1)$ and $I\cap R_1$ is a subring of the ring 
$(\,R_1, \, +, \, \,\sharp\, \,)$.
\end{description}
\end{definition}

\medskip
Let $I$ be a triideal of a triring $(\,R=R_0\oplus R_1, \, +, \,\cdot, \, \,\sharp\, \,)$. It is clear that $\displaystyle\frac{R}{I}=\left(\displaystyle\frac{R}{I}\right)_0\oplus 
\left(\displaystyle\frac{R}{I}\right)_1$  with
$\left(\displaystyle\frac{R}{I}\right)_i:=\displaystyle\frac{R_i+I}{I}$ for $i=0$ and $1$. We now define a local product on $\left(\displaystyle\frac{R}{I}\right)_1$ by
\begin{equation}\label{eq5}
(\alpha +I)\,\sharp\, (\beta +I):=\alpha\,\sharp\,\beta +I 
\quad\mbox{for $\alpha$, $\beta\in R_1$.}
\end{equation}
Then the local product defined by (\ref{eq5}) is well-defined, and the triassociative law holds. 
Therefore, $\displaystyle\frac{R}{I}$ becomes a triring, which is called the {\bf quotient triring} of $R$ with respect to the triideal $I$.

\medskip
\begin{definition}\label{def1.4} Let $R=R_0\oplus R_1$ and 
$\overline{R}=\overline{R}_0\oplus \overline{R}_1$ be trirings. A map 
$\phi: R\to \overline{R}$ is called a {\bf triring homomorphism} if 
$$
\phi (x+y)=\phi (x)+\phi (y), \quad \phi (xy)=\phi (x)\phi (y), \quad
\phi (1_R)=1_{\overline{R}},
$$$$
\phi (R_0)\subseteq \overline{R}_0, \quad \phi (R_1)\subseteq \overline{R}_1,
$$$$
\phi (\alpha \,\sharp\, \beta )=\phi (\alpha) \,\sharp\, \phi (\beta), \quad
\phi(1^{\sharp})=\bar{1}^{\sharp}
$$
where $x$, $y\in R$, $\alpha$, $\beta\in R_1$, $1_R$ and $1_{\overline{R}}$ are the identities of $R$ and $\overline{R}$ respectively, and $1^{\sharp}$ and $\bar{1}^{\sharp}$ are the local identities of $R$ and $\overline{R}$ respectively. A bijective triring homomorphism is called a {\bf triring isomorphism}. We shall use $R\simeq \overline{R}$ to indicate that there exists a triring isomorphism from $R$ to $\overline{R}$.
\end{definition}

\medskip
Let  $\phi$ be a  triring homomorphism from a triring $R$ to a triring $\overline{R}$. the {\bf kernel} $Ker\phi$ and the {\bf image} $Im\phi$ of $\phi$ are defined by
$$ Ker\phi :=\{\, x \,|\, \mbox{$a\in R$ and $\phi (x)=0$} \, \}$$ 
and
$$ Im\phi :=\{\, \phi (x) \,|\, \mbox{$x\in R$} \, \}.$$ 
Clearly, $Ker\phi$ is a triideal of $R$, $Im\phi$ is a subtriring of $\overline{R}$ and 
$$\overline{\phi}: x+Ker\phi \to \phi (x) \quad\mbox{for $x\in R$}$$
is a triring isomorphism from the quotient triring $\displaystyle\frac{R}{Ker\phi}$ to the subtriring $Im\phi$ of $\overline{R}$. If $I$ is triideal of $R$, then 
$$ \nu : x\mapsto x+I \quad\mbox{for $x\in R$}$$
is a surjective triring homomorphism from $R$ the quotient tiring $\displaystyle\frac{R}{I}$ with kernel $I$. The map $\nu$ is called the {\bf natural triring homomorphism}.

\bigskip
\begin{proposition}\label{pr1.1}  Let $\phi$ be a surjective homomorphism from a triring $R=R_0\oplus R_1$ to a triring $\overline{R}=\overline{R}_0\oplus \overline{R}_1$.
\begin{description}
\item[(i)] $\phi (R_i)=\overline{R}_i$ for $i=0$ and $1$.
\item[(ii)] Let
$$
\mathcal{S}: =\{\, I \,|\, \mbox{$I$ is a triideal of $R$ and $I \supseteq Ker \phi $} \, \}
$$
and
$$
\overline{\mathcal{S}}: =\{\, \overline{I} \,|\, \mbox{$\overline{I}$ is a triideal of $\overline{R}$} \, \}.
$$
The map 
$$
\Psi : I \mapsto \phi (I): =\{\, \phi (x) \,|\, x\in I \, \}
$$
is a bijection from $\mathcal{S}$ to $\overline{\mathcal{S}}$, and the inverse map $\Psi ^{-1} : \overline{\mathcal{S}} \to \mathcal{S}$ is given by
$$
\Psi ^{-1} : \overline{I} \mapsto \phi ^{-1} (\overline{I}) : = \{\, x \,|\, \mbox{$x\in R$ and $\phi (x)\in \overline{I}$ \, \}.}
$$
\item[(iii)] If $I$ is a triideal of $R$ containing $Ker \phi$, then the map
$$
x+I \mapsto \phi (x) + \phi (I) \quad\mbox{for $x\in R$}
$$
is a triring isomorphism from the quotient triring $\displaystyle\frac{R}{I}$ onto the quotient triring $\displaystyle\frac{\overline{R}}{\phi (I)}$.
\end{description}
\end{proposition}

\medskip
\noindent
{\bf Proof} A routine check.

\hfill\raisebox{1mm}{\framebox[2mm]{}}

\medskip
We now prove a basic property for trirings.

\begin{proposition}\label{pr2.1}  Let $(\,R=R_0\oplus R_1, \, +, \,\cdot, \, \,\sharp\, \,)$ be a triring. If $x_i\in R_i$  for $i\in\{\, 0, \, 1\,\}$, then both $Rx_0=R_0x_0\oplus R_1x_0$ and $R_1\,\sharp\, x_1$ are triideals of $R$.
\end{proposition}

\medskip
\noindent
{\bf Proof} Since $(R_1, \, +, \, \sharp)$ is a commutative ring, $R_1\,\sharp\, x_1$ is clearly a triideal of $R$.

\medskip
Using the triassociative law and (\ref{eq6}), we have
\begin{equation}\label{eq7}
\Big(R_i(R_0x_0+R_1x_0)\Big)\bigcup
\Big(R_0x_0+R_1x_0)R_i\Big)\subseteq R_0x_0+R_1x_0
\end{equation}
for $i=0$, $1$. By (\ref{eq7}), $R_0x_0\oplus R_1x_0$ is an ideal of the ring $(\,R, \, +, \, \cdot\,)$. Also, $R_1x_0=R_1\,\sharp\, (1^{\sharp}x_0)$ is obviously an ideal of the commutative ring $(\,R_1, \, +, \, \sharp\,)$. Thus, $R_0x_0\oplus R_1x_0$ is a triideal of $R$.

\hfill\raisebox{1mm}{\framebox[2mm]{}}

\medskip
The next position gives some operations about triideals in a triring.

\medskip
\begin{proposition}\label{pr2.2}  Let $I$, $J$, $I_{\lambda}$ with $\lambda\in \Lambda$ be triideals of a triring $R$.
\begin{description}
\item[(i)] The intersection $I\cap J$ and the sum 
$\displaystyle\sum_{\lambda\in \Lambda}I_{\lambda}$ are triideals of $R$.Moreover, we have
$(I\cap J)_i=I_i\cap J_i$ and 
$\left(\displaystyle\sum_{\lambda\in \Lambda}I_{\lambda}\right)_i=
\displaystyle\sum_{\lambda\in \Lambda}(I_{\lambda})_i$ for $i=0$, $1$.
\item[(ii)] The {\bf mixed product} 
$I\,\stackrel{\cdot}{\sharp}\,J:=(I\,\stackrel{\cdot}{\sharp}\,J)_0\oplus 
(I\,\stackrel{\cdot}{\sharp}\,J)_1$ of $I$ and $J$ is a triideal, where 
$(I\,\stackrel{\cdot}{\sharp}\,J)_0=I_0J_0$ and $(I\,\stackrel{\cdot}{\sharp}\,J)_1=I_1\,\sharp\,J_1$.
\end{description}
\end{proposition}

\medskip
\noindent
{\bf Proof} Clear.

\hfill\raisebox{1mm}{\framebox[2mm]{}}

\bigskip
\section{Prime Triideals}

\medskip
We begin this section by introducing the notion of prime triideals.

\begin{definition}\label{def5.1} Let $(\, R=R_0\oplus R_1, \, + , \, \cdot\, , \,\sharp\, \,)$ be a triring. An triideal $P=P_0\oplus P_1$ of $R$ is called a {\bf prime triideal} if $P\ne R$ and 
\begin{eqnarray}
\label{eq15}x_0y_0\in P_0  &\Rightarrow&  \mbox{$x_0\in P_0$ or $y_0\in P_0$},\\
\label{eq16}x_0y_1\in P_1  &\Rightarrow&  \mbox{$x_0\in P_0$ or $y_1\in P_1$},\\
\label{eq17}x_1y_0\in P_1  &\Rightarrow&  \mbox{$x_1\in P_1$ or $y_0\in P_0$},\\
\label{eq18}x_1\,\sharp\, y_1\in P_1 &\Rightarrow&  \mbox{$x_1\in P_1$ or $y_1\in P_1$},
\end{eqnarray}
where $x_i$, $y_i\in R_i$ for $i=0$ and $1$.
\end{definition}

\medskip
Let $(\, R=R_0\oplus R_1, \, + , \, \cdot\, \,\sharp\, \,)$ be a triring. The set of all prime triideals of $R$ is called the {\bf trispectrum} of $R$ and denoted by $Spec^{\sharp}R$. It is clear that
$$Spec^{\sharp}R=Spec^{\sharp}_0R\cup Spec^{\sharp}_1R\quad\mbox{and}\quad 
Spec^{\sharp}_0R\cap Spec^{\sharp}_1R=\emptyset,$$
where
$$Spec^{\sharp}_0R:=\{\, P \, |\, \mbox{$P\in Spec^{\sharp}R$ and $P\supseteq R_1$}\,\}$$
is called the {\bf even trispectrum} of $R$ and 
$$Spec^{\sharp}_1R:=\{\, P \, |\, \mbox{$P\in Spec^{\sharp}R$ and $P\not\supseteq R_1$}\,\}$$
is called the {\bf odd trispectrum} of $R$.  

\medskip
Let $P=P_0\oplus P_1$ be a triideal of a triring $R$. Then 
$P\in Spec^{\sharp}_0R$ if and only if $P_0$ is a prime ideal of the commutative ring 
$(\, R_0, \, + , \, \cdot\, \,)$. It is also obvious that $P\in Spec^{\sharp}_1R$ if and only if
$P_0$ is a prime ideal of the commutative ring $(\, R_0, \, + , \, \cdot\, \,)$, $P_1$ is a prime ideal of the commutative ring $(\, R_1 , \, + , \, \,\sharp\, \, )$, and the $\left(\displaystyle\frac{R_0}{P_0}, \displaystyle\frac{R_0}{P_0}\right)$-bimodule $\displaystyle\frac{R}{P}$ is faithful as both left module and right module, where the left $\displaystyle\frac{R_0}{P_0}$-module action on $\displaystyle\frac{R}{P}$ is defined by 
$$
(x_0+P_0)(y+P):= x_0y+P \quad\mbox{for $x_0\in P_0$ and $y\in R$}
$$
and the right $\displaystyle\frac{R_0}{P_0}$-module action  on $\displaystyle\frac{R}{P}$ is defined by 
$$
(y+P)(x_0+P_0):= yx_0+P \quad\mbox{for $x_0\in P_0$ and $y\in R$}.
$$

\medskip
Clearly, the even trispectrum $Spec^{\sharp}_0R$ of a triring $R$ is not empty. A basic property of trirings is that the odd trispectrum $Spec^{\sharp}_1R$ is always not empty provided $R_1\ne 0$. This basic fact is a corollary of the following

\begin{proposition}\label{pr5.1} Let $(\, R=R_0\oplus R_1, \, + , \, \cdot\, \,\sharp\, \,)$ be a triring with $R_1\ne 0$. If $P_1$ is a prime ideal of the commutative ring 
$(\, R_1 , \, + , \, \,\sharp\, \, )$, then there exists prime ideal $P_0$ of the commutative ring $(\, R_0 , \, + , \, \,\cdot\, \, )$ such that $P:=P_0\oplus P_1$ is a prime triideal of $R$, and $P_0$ contains every ideal $I_0$ of the ring  $(\, R_0 , \, + , \, \,\cdot\, \, )$ which has the property: $R_1I_0\subseteq P_1$.
\end{proposition}

\medskip
\noindent
{\bf Proof} Consider the set $\Omega$ defined by
$$
\Omega :=\{\,I_0\,|\, \mbox{$I_0$ is an ideal of $(\, R_0 , \, + , \, \,\cdot\, \, )$ and $R_1I_0\subseteq P_1$}\,\}.
$$

Clearly, $1\not\in I_0$ if $I_0\in \Omega$. Since $0\in \Omega$, $\Omega$ is not empty. The relation of inclusion, $\subseteq$, is a partial order on $\Omega$. Let $\Delta$ be a non-empty totally ordered subset of $\Omega$. Let $J_0:=\displaystyle\bigcup _{I_0\in\Delta} I_0$. Then $J_0\in\Omega$. Thus $J_0$ is an upper bound for $\Delta$ in $\Omega$. By Zorn's Lemma, the partial order set $(\Omega,\,\subseteq)$ has a maximal element $P_0$. We are going to prove that $P:=P_0\oplus P_1$ is a prime triideal of $R$. 
Clearly, $P=P_0\oplus P_1$ is a triideal satisfying (\ref{eq18}). Let $x_i$, $y_i\in R_i$ for $i=0$, $1$. 

\medskip
If $x_1y_0\in P_1$ and $x_1\not\in P_1$, then 
$x_1\,\sharp\,(1^{\sharp}y_0)=x_1y_0\in P_1$, which implies that $1^{\sharp}y_0\in P_1$. Hence, we have
\begin{eqnarray}\label{eq19}
R_1(P_0+R_0y_0)&\subseteq& R_1P_0+R_1R_0y_0\subseteq P_1+R_1y_0\nonumber\\
&=&P_1+R_1\,\sharp\,(1^{\sharp}y_0)\subseteq P_1+R_1\,\sharp\,P_1\subseteq P_1.
\end{eqnarray}
Using (\ref{eq19}) and the fact that $P_0+R_0y_0$ is an ideal of $R_0$, we get
$P_0+R_0y_0\in \Omega$. Since $P_0+R_0y_0\supseteq P_0$, we have to have 
$P_0+R_0y_0= P_0$, which implies that $y_0\in P_0$. This proves that (\ref{eq17}) holds.
Similarly, we have that both (\ref{eq15}) and (\ref{eq16}) hold. Hence, $P=P_0\oplus P_1$ is a prime triideal of $R$. 

\hfill\raisebox{1mm}{\framebox[2mm]{}}

\medskip
The following proposition gives another characterization of prime triideals.

\begin{proposition}\label{pr5.2} Let $(\, R=R_0\oplus R_1, \, + , \, \cdot , \, \,\sharp\, \,)$ be a triring. The following are equivalent.
\begin{description}
\item[(i)] $P$ is a prime triideal.
\item[(ii)] For two triideals $I$, $J$ of $R$, $I\,\stackrel{\cdot}{\sharp}\,J\subseteq P$ implies that $I\subseteq P$ or $J\subseteq P$.
\end{description}
\end{proposition}

\medskip
\noindent
{\bf Proof} Use By Proposition \ref{pr2.1}.

\hfill\raisebox{1mm}{\framebox[2mm]{}}

\bigskip
\section{Trinilradicals}

Let $R=R_0\oplus R_1$ be a triring with the local identity $1^\sharp$. For 
$\alpha \in R_1$, the {\bf local $n$th power} $\alpha ^{\sharp n}$ is defined by:
$$
\alpha^{\sharp n} :=\left\{ \begin{array}{ll}
1^\sharp, & \qquad\mbox{if $n=0$;}\\
\underbrace{\alpha \,\sharp\, \alpha  \,\sharp\, \cdots  \,\sharp\, \alpha}_n , & 
\qquad\mbox{if $n$ is a positive integer.} \end{array}\right.
$$
The products $(x^m)(\alpha^{\sharp n})$ and $(\alpha^{\sharp n})(x^m)$ will be denoted by $x^m\alpha^{\sharp n}$ and $\alpha^{\sharp n}x^m$ respectively, where $x\in R$ and $\alpha \in R_1$.

\begin{proposition}\label{pr6.1} Let $(\, R=R_0\oplus R_1, \, + , \, \cdot , \, \,\sharp\, \,)$ be a triring with the local identity $1^\sharp$. 
If $x$, $y\in R$, $\alpha$, $\beta\in R_1$ and $m\in\mathcal{Z}_{>0}$, then
\begin{equation}\label{eq27}
(x\alpha)\,\sharp\, (y\beta)=(xy)(\alpha\,\sharp\, \beta), \qquad
(\alpha x)\,\sharp\, (\beta y)=(\alpha\,\sharp\, \beta)xy
\end{equation}
and
\begin{equation}\label{eq28}
(x\alpha)^{\sharp m}=x^m \alpha ^{\sharp m}, \qquad 
(\alpha x)^{\sharp m}=\alpha ^{\sharp m}x^m.
\end{equation}
\end{proposition}

\medskip
\noindent
{\bf Proof} By the triassociative law, we have
$$
(x\alpha )\,\sharp\, (y\beta)=x\left(\alpha\,\sharp\, (y\beta)\right)=
x\left((y\beta)\,\sharp\, \alpha\right)=(xy)(\beta\,\sharp\, \alpha)=(xy)(\alpha\,\sharp\, \beta)
$$
and
$$
(\alpha x)\,\sharp\, (\beta y)=\left((\alpha x)\,\sharp\, \beta)\right)y=
\left((\beta)\,\sharp\, (\alpha x)\right)y=(\beta\,\sharp\, \alpha)xy=
(\alpha\,\sharp\, \beta)xy.
$$
Hence, (\ref{eq27}) holds. Clearly, (\ref{eq28}) follows from (\ref{eq27}).

\hfill\raisebox{1mm}{\framebox[2mm]{}}

\medskip
\begin{definition}\label{def6.1} Let $(\, R=R_0\oplus R_1, \, + , \, \cdot , \, \,\sharp\, \,)$ be a triring. 
\begin{description}
\item[(i)] An element $x$ of $R$ is said to be {\bf  trinilpotent} if 
\begin{equation}\label{eq29}
\mbox{$x_0^m=0$ and $x_1^{\,\sharp\, n}=0$ for some $m$, $n\in\mathcal{Z}_{>0}$,}
\end{equation}
where $x_0$ and $x_1$ are the even component and the old component of $x$ respectively.
\item[(ii)] The set of all trinilpotent elements of $R$ is called the {\bf trinilradical} of $R$ and denoted by $nilrad ^{\sharp}(R)$ or $\sqrt[\sharp]{0}$.
\end{description}
\end{definition}

\medskip
The ordinary nilradical of a ring $(\, A, \, + , \, \cdot  \,)$ is denoted by $nilrad(A)$ or
$nilrad(A, + , \cdot)$; that is,
$$
nilrad(A):=\{\, x \,|\, \mbox{$x^m=0$ for some $m\in\mathcal{Z}_{>0}$}\,\}.
$$

If $(\, R=R_0\oplus R_1, \, + , \, \cdot , \, \,\sharp\, \,)$ is a triring, then
$$
nilrad(R, +, \cdot)=nilrad(R_0, + , \cdot)\oplus R_1
$$
and
\begin{equation}\label{eq30}
nilrad^{\sharp}R=nilrad(R_0, + , \cdot)\oplus nilrad(R_1, +, \,\sharp\,).
\end{equation}
Hence, $nilrad(R, +, \cdot)\supseteq nilrad^{\sharp}R$, i.e., the trinilradical of a triring $R$ is smaller than the ordinary nilradical of the ring $(R,\,+, \,\cdot)$.

\begin{proposition}\label{pr6.2} Let $(\, R=R_0\oplus R_1, \, + , \, \cdot , \, \,\sharp\, \,)$ be a triring.
\begin{description}
\item[(i)] The trinilradical $nilrad ^{\sharp}(R)$ is a triideal of $R$.
\item[(ii)] $nilrad ^{\sharp}\left( \displaystyle\frac{R}{nilrad ^{\sharp}(R)}\right)=0.$
\end{description}
\end{proposition}

\medskip
\noindent
{\bf Proof} (i) By the definition of trinilradicals, we have 
$$
(nilrad^{\sharp}R)\cap R_0=nilrad(R_0, + , \cdot)\quad\mbox{and}\quad
(nilrad^{\sharp}R)\cap R_1=nilrad(R_1, +, \,\sharp\,).
$$
Using the fact above and (\ref{eq30}), we need only to prove
\begin{equation}\label{eq31}
\{xa,\, ax\}\subseteq nilrad ^{\sharp}R\quad\mbox{for $x\in R$ and $a\in nilrad ^{\sharp}R$}.
\end{equation}

Let $x=x_0+x_1\in R$ and $a=a_0+a_1\in nilrad ^{\sharp}R$, where $x_i$, $a_i\in R_i$ for $i=0$, $1$. Then we have $a_0^m=a_1^{\,\sharp\, m}=0$ for some $m\in\mathcal{Z}_{>0}$. Since
\begin{eqnarray}
\label{eq32} xa&=&(x_0+x_1)(a_0+a_1)=x_0a_0+(x_0a_1+x_1a_0),\\
(x_0a_0)^m&=&x_0^ma_0^m=x_0^m0=0,\nonumber\\
(x_0a_1)^{\sharp m}&=&x_0^ma_1^{\sharp m}=x_0^m0=0,\nonumber\\
(x_1a_0)^{\sharp m}&=&(x_1a_0)^{\sharp m}=x_1^{\sharp m}a_0^m=x_1^{\sharp m}0=0,\nonumber
\end{eqnarray}
we get
\begin{equation}\label{eq33}
x_0a_0\in nilrad(R_0, \,+, \,\cdot) \quad\mbox{and}\quad 
x_0a_1+x_1a_0\in nilrad ^{\sharp}(R_1, \,+,\, \sharp).
\end{equation}

It follows from (\ref{eq32}) and (\ref{eq33}) that $xa\in nilrad ^{\sharp}R$. Similarly, we have $ax\in nilrad ^{\sharp}R$. This proves (i).

\medskip
(ii) follows from (i).

\hfill\raisebox{1mm}{\framebox[2mm]{}}

\medskip
If $I$ is a triideal of a triring $R=R_0\oplus R_1$, then the {\bf trinilradical} $\sqrt[\sharp]{I}$ of $I$ is defined by
$$
\sqrt[\sharp]{I}:=\{\, x\in R \, | \, \mbox{$x_0^m\in I_0$ and $x_1^{\,\sharp\, n}\in I_1$ for some $m$, $n\in\mathcal{Z}_{>0}$} \,\},
$$
where $x=x_0+x_1$, $x_i\in R_i$ and $I_i=I\cap R_i$ for $i=0$, $1$. Since
$$
nilrad ^{\sharp}\left(\displaystyle\frac{R}{I}\right)=\displaystyle\frac{\sqrt[\sharp]{I}}{I},
$$
$\sqrt[\sharp]{I}$ is a triideal of $R$. A triideal $I$ of a triring $R$ is called a {\bf radical triideal} if $\sqrt[\sharp]{I}=I$.

\medskip
We now characterize the trinilradical of a triring by using prime triideals. 

\begin{proposition}\label{pr6.3} Let $(\, R=R_0\oplus R_1, \, + , \, \cdot\, , \,\sharp\, \,)$ be a triring. The trinilradical of $R$ is the intersection of the prime triideals of $R$.
\end{proposition}

\medskip
\noindent
{\bf Proof}  Let $x=x_0+x_1$ be any element of $nilrad ^{\sharp}(R)$, where $x_0\in R_0$ and $x_1\in R_1$. Then $x_0^m=0$ and $x_1^{\,\sharp\, n}=0$ for some $m$, $n\in\mathcal{Z}_{>0}$. Let $P=P_0\oplus P_1$ be any prime triideal of $R$. Since $x_0^m=0\in P_0$, we have 
\begin{equation}\label{eq34} 
x_0\in P_0
\end{equation}
by (\ref{eq15}).

If $P\not\supseteq R_1$, then $P_1$ is a prime ideal of the commutative ring 
$(\, R_1, \, + , \, \,\sharp\, \,)$. Using this fact and $x_1^{\,\sharp\, n}=0$, we get
\begin{equation}\label{eq35} 
x_1\in P_1.
\end{equation}
If $P\supseteq R_1$, then (\ref{eq35}) is obviously true. It follows from (\ref{eq34}) and (\ref{eq35}) that $x=x_0+x_1\in P$. This proves that
\begin{equation}\label{eq36} 
nilrad ^{\sharp}(R)\subseteq \bigcap_{P\in Spec^{\sharp}R} P.
\end{equation}

\medskip
Conversely, we prove that
\begin{equation}\label{eq37} 
z\not\in nilrad ^{\sharp}(R) \,\Rightarrow\, z\not\in \bigcap_{P\in Spec^{\sharp}R} P.
\end{equation}

\medskip
\underline{\it Case 1:} $z^m\ne 0$ for all $m\in \mathcal{Z}_{>0}$, in which case, 
$z^m\not\in R_1$ for all $m\in \mathcal{Z}_{>0}$. Hence, $z+ R_1$ is not a nilpotent element of the commutative ring $\displaystyle\frac{R}{R_1}$; that is,
$$
z+ R_1\ne nilrad\left(\displaystyle\frac{R}{R_1} \right)
=\bigcap _{\frac{I}{R_1}\in Spec \left( \frac{R}{R_1}\right)} 
\left(\displaystyle\frac{I}{R_1}\right),
$$
where $Spec \left( \displaystyle\frac{R}{R_1}\right)$ is the ordinary spectrum of the commutative ring 
$\displaystyle\frac{R}{R_1}$. Hence, there exists a prime ideal $\displaystyle\frac{I}{R_1}$ of the commutative ring $\displaystyle\frac{R}{R_1}$ such that $z\not\in I$. Since $I$ is a prime triideal of $R$, (\ref{eq37}) holds in this case.

\medskip
\underline{\it Case 2:} $z^m= 0$ for some $m\in \mathcal{Z}_{>0}$. Let $z=z_0+z_1$ with 
$z_0\in R_0$ and $z_1\in R_1$. Then $0=z^m=z_0^m+z_0^{m-1}r_1$ for some $r_1\in R_1$ by (\ref{eq6}). Thus, $z_0^m=0$, which implies that  $z_1^{\,\sharp\, n}\ne 0$ for all 
$n\in \mathcal{Z}_{>0}$ in this case.
We  now consider the following set
$$
T:=\left\{\, J \, \left| \, \mbox{$J$ is a triideal of $R$ and $z_1^{\,\sharp\, n}\not\in J$ for all $n\in \mathcal{Z}_{>0}$}\right.\right\}.
$$
Since $\{0\}\in T$, $T$ is nonempty. Clearly, $(\,T, \, \subseteq \,)$ is a partially order set, where $\subseteq $ is the relation of set inclusion. If $\{\, J_{\lambda} \, |\, \lambda\in \Lambda \,\}$ is a nonempty totally ordered subset of $T$, then 
$\bigcup_{\lambda\in \Lambda} J_{\lambda}$ 
is an upper bound of 
$\{\, J_{\lambda} \, |\, \lambda\in \Lambda \,\}$ in $T$. By Zorn's Lemma, the partially ordered set $(\,T, \, \subseteq \,)$ has a maximal element $P$. We are going to prove that $P$ is a prime triideal of $R$.

\medskip
Let $x=x_0+x_1$ and $y=y_0+y_1$ be two elements of $R$, where $x_i$, $y_i\in R_i$ for $i=0$, $1$. First, if $x_0\not\in P_0$ and $y_0\not\in P_0$, then 
\begin{equation}\label{eq38}
P\subset P+Rx_0\quad\mbox{and}\quad P\subset P+Ry_0.
\end{equation}
Since both $P+Rx_0$ and $P+Ry_0$ are triideals of $R$ by Proposition~\ref{pr2.1}, (\ref{eq38}) implies that
$$
z_1^{\,\sharp\, u}\in P+Rx_0\quad\mbox{and}\quad z_1^{\,\sharp\, v}\in P+Ry_0
$$
or
$$
z_1^{\,\sharp\, u}\in P_1+R_1x_0\quad\mbox{and}\quad z_1^{\,\sharp\, v}\in P_1+R_1y_0\quad\mbox{for some $u$, $v\in \mathcal{Z}_{>0}$.}
$$
Thus, we have
\begin{eqnarray*}
&&z_1^{\,\sharp\, (u+v)}=z_1^{\,\sharp\, u}\,\sharp\, z_1^{\,\sharp\, v}\in (P_1+R_1x_0)\,\sharp\, (P_1+R_1y_0)\\&&\\
&\subseteq& 
\underbrace{P_1\,\sharp\, P_1+ P_1\,\sharp\, (R_1y_0)+
(R_1x_0)\,\sharp\, P_1}_{\mbox{This is a subset of $P$}}+ (R_1x_0)\,\sharp\, (R_1y_0)\\&&\\
&\subseteq& P+R_1x_0y_0,
\end{eqnarray*}
which implies that $x_0y_0\not\in P$. This proves that
\begin{equation}\label{eq39} 
\mbox{$x_0\not\in P_0$ and $y_0\not\in P_0$} \,\Rightarrow\, x_0y_0\not\in P_0.
\end{equation}

Similarly, we have
\begin{equation}\label{eq44} 
\mbox{$x_0\not\in P_0$ and $y_1\not\in P_1$} \,\Rightarrow\, x_0y_1\not\in P_1,
\end{equation}

\begin{equation}\label{eq47} 
\mbox{$y_1\not\in P_1$ and $x_0\not\in P_0$} \,\Rightarrow\, y_1x_0\not\in P_1
\end{equation}
and
\begin{equation}\label{eq48} 
\mbox{$x_1\not\in P_1$ and $y_1\not\in P_1$} \,\Rightarrow\, x_1y_1\not\in P_0.
\end{equation}

\medskip
By (\ref{eq39}), (\ref{eq44}), (\ref{eq47}) and (\ref{eq48}), $P$ is a prime triideal. Since $z_1\not\in P$, (\ref{eq37}) also holds in Case 2.

\medskip
It follows from (\ref{eq36}) and (\ref{eq37}) that Proposition~\ref{pr6.3} is true.

\hfill\raisebox{1mm}{\framebox[2mm]{}}

\bigskip
The next proposition is a corollary of Proposition~\ref{pr6.3}.

\begin{proposition}\label{pr6.4} If $I$ is a triideal of a Hu-Liu triring $R$ and $I\ne R$, then
$$\sqrt[\sharp]{I}=\bigcap _{\mbox{$P\in Spec^{\sharp}R$ and $P\supseteq I$}} P.$$
\end{proposition}

\medskip
\noindent
{\bf Proof} By Proposition~\ref{pr6.3}, we have
\begin{eqnarray*}
x\in \sqrt[\sharp]{I}
&\Leftrightarrow& x+I\in nilrad ^{\sharp}\left(\frac{R}{I}\right)
=\bigcap _{\mbox{$\frac{P}{I}\in Spec^{\sharp}\left(\frac{R}{I}\right)$}}\frac{P}{I}\\
&\Leftrightarrow& x\in \bigcap _{\mbox{$P\in Spec^{\sharp}R$ and $P\supseteq I$}} P.
\end{eqnarray*}

\hfill\raisebox{1mm}{\framebox[2mm]{}}

\bigskip
\section{Extended Zariski Topology}

Let $(\, R=R_0\oplus R_1, \, + , \, \cdot , \, \,\sharp\, \,)$ be a triring. For a triideal $I$ of $R$, we define a subset $\mathcal{V}(I)$ of $Spec^{\sharp}R$ by
\begin{equation}\label{eq49} 
V^{\sharp}(I): =\{\, P \, |\, \mbox{$P\in Spec^{\sharp}R$ and $P\supseteq I$} \,\}.
\end{equation}

\begin{proposition}\label{pr7.1} Let $R$ be a triring.
\begin{description}
\item[(i)] $V^{\sharp}(0)=spec^{\sharp}R$ and $V^{\sharp}(R)=\emptyset$.
\item[(ii)] $V^{\sharp}(I)\cup V^{\sharp}(J)
=V^{\sharp}(I\cap J)=V^{\sharp}(I\,\stackrel{\cdot}{\sharp}\,J)$, where $I$ and $J$ are two triideals of $R$.
\item[(iii)] $\displaystyle\bigcap _{\lambda \in \Lambda} 
V^{\sharp}\left(I_{(\lambda)}\right)=
V^{\sharp}\left(\displaystyle\sum _{\lambda \in \Lambda}I_{(\lambda)}\right)$, where 
$\left\{\left.\,I_{(\lambda)}\, \right| \, \lambda \in \Lambda\, \right\}$ is a set of triideals of $R$.
\end{description}
\end{proposition}

\medskip
\noindent
{\bf Proof} Since (i) and (iii) are clear, we need only to prove (ii).

\medskip
Since $I\,\stackrel{\cdot}{\sharp}\,J\subseteq I\cap J\subseteq I$, we get $V^{\sharp}(I\,\stackrel{\cdot}{\sharp}\,J)\supseteq V^{\sharp}(I\cap J)
\supseteq V^{\sharp}(I)$. Similarly, we have $V^{\sharp}(I\,\stackrel{\cdot}{\sharp}\,J)\supseteq V^{\sharp}(I\cap J)\supseteq V^{\sharp}(J)$. Thus, we get
\begin{equation}\label{eq50} 
V^{\sharp}(I)\cup V^{\sharp}(J)\subseteq V^{\sharp}(I\cap J)
\subseteq V^{\sharp}(I\,\stackrel{\cdot}{\sharp}\,J).
\end{equation}

\medskip
Conversely, if $P\in V^{\sharp}(I\,\stackrel{\cdot}{\sharp}\,J)$, then $I\,\stackrel{\cdot}{\sharp}\,J\subseteq P$. By Proposition~\ref{pr5.2}, we get $I\subseteq P$ or 
$J\subseteq P$. Hence, $P\in V^{\sharp}(I)\cup V^{\sharp}(J)$. This proves that
\begin{equation}\label{eq51} 
V^{\sharp}(I\,\stackrel{\cdot}{\sharp}\,J)\subseteq  V^{\sharp}(I)\cup V^{\sharp}(J).
\end{equation}

\medskip
It follows from (\ref{eq50}) and (\ref{eq51}) that (ii) is true.

\hfill\raisebox{1mm}{\framebox[2mm]{}}

\bigskip
Let $(\, R=R_0\oplus R_1, \, + , \, \cdot\, , \,\sharp\, \,)$ be a triring. By Proposition~\ref{pr7.1}, the collection
$$
V^{\sharp}: =\{\, V^{\sharp}(I) \, |\, \mbox{$I$ is a triideal of $R$} \, \}
$$
of subsets of $Spec^{\sharp}R$ satisfies the axioms for closed sets in a topological space. The topology on $Spec^{\sharp}R$ having the elements of $V^{\sharp}$ as closed sets is called the {\bf extended Zariski topology}. The collection
$$
D^{\sharp}: =\{\, D^{\sharp}(I) \, |\, \mbox{$I$ is a triideal of $R$} \, \}
$$
consists of the open sets of the extended Zariski topology on $Spec^{\sharp}R$, where
$$
D^{\sharp}(I): =Spec^{\sharp}R\setminus V^{\sharp}(I)=
\{\, P \, |\, \mbox{$P\in Spec^{\sharp}R$ and $P\not\supseteq I$} \, \}.
$$

\medskip
For $x_i\in R_i$ with $i\in \{0,\, 1\}$, both $Rx_0$ and $R_1\,\sharp\,x_1$ are triideals of $R$ by Proposition~\ref{pr2.1}. Let
$$
D^{\sharp}(x_0): =D^{\sharp}(Rx_0), \qquad D^{\sharp}(x_1): =D^{\sharp}(R_1\,\sharp\,x_1).
$$
If $I_0$ and $I_1$ are the even part and odd part of an triideal $I$, then 
$$D^{\sharp}(I)=\bigcup _{x_i\in I_i,\, i=0,1}D^{\sharp}(x_i).$$
Thus, $\{\, D^{\sharp}(x_i)\,|\,\mbox{$x_i\in R_i$ with  $i=0,\, 1$}\,\}$ forms an open base for the extended Zariski topology on $Spec^{\sharp}R$. Each $D^{\sharp}(x_i)$ is called a {\bf basic open subset} of $Spec^{\sharp}R$. Clearly, $D^{\sharp}(0)=\emptyset$, $D^{\sharp}(1)=Spec^{\sharp}R$, $D^{\sharp}(1^{\sharp})=Spec^{\sharp}_1R$, and $D^{\sharp}(x_1)\subseteq  Spec^{\sharp}_1R$ for $x_1\in R_1$.

\begin{proposition}\label{pr7.2} Let $I$ and $J$ be triideals of a Hu-Liu triring. 
\begin{description}
\item[(i)] $V^{\sharp}(I)\subseteq V^{\sharp}(J)$ if and only if 
$\sqrt[\sharp]{J}\subseteq \sqrt[\sharp]{I}$.
\item[(ii)] $V^{\sharp}(I)=V^{\sharp}(\sqrt[\sharp]{I})$.
\end{description}
\end{proposition}

\medskip
\noindent
{\bf Proof} (i) follows from Proposition~\ref{pr6.4}, and (ii) follows from (i).

\hfill\raisebox{1mm}{\framebox[2mm]{}}

\bigskip
\begin{definition}\label{def7.1} Let $X$ be a topological space.
\begin{description}
\item[(i)] A closed subset $F$ of $X$ is {\bf reducible} if $F=F_{(1)}\cup F_{(2)}$ for proper closed subsets $F_{(1)}$, $F_{(2)}$ of $X$. We call a closed subset $F$ {\bf irreducible} if it is not reducible.
\item[(ii)] $X$ is {\bf quasicompact} if given an arbitrary open covering 
$\{\, U_{(i)}\,|\, i\in I\,\}$ of $X$, there exists a finite subcovering of $X$, i.e., there exist finitely many members $U_{(i_1)}$, $\dots$, $U_{(i_n)}$ of $\{\, U_{(i)}\,|\, i\in I\,\}$ such that $X=U_{(i_1)}\cup\cdots\cup U_{(i_n)}$.
\end{description}
\end{definition}

\medskip
Using the topological concepts above, we have the following

\medskip
\begin{proposition}\label{pr7.3}  Let $(\, R=R_0\oplus R_1, \, + , \, \cdot , \, \,\sharp\, \,)$ be a triring.
\begin{description}
\item[(i)] The trispectrum $Spec^{\sharp}R$ is quasicompact.
\item[(ii)] Both $Spec^{\sharp}_0R$ and $Spec^{\sharp}_1R$ are quasicompact subsets of $Spec^{\sharp}R$.
\item[(iii)] If $I$ is a triideal of $R$, then the closed subset $V^{\sharp}(I)$ of $Spec^{\sharp}R$ is irreducible if and only if $\sqrt[\sharp]{I}$ is a prime triideal.
\end{description}
\end{proposition}

\medskip
\noindent
{\bf Proof}  (i) Let $\{\, D^{\sharp}(I_{(i)})\,|\, i\in\Delta\,\}$ be an open covering of $Spec^{\sharp}R$, where $I_{(i)}=I_{(i)0}\oplus I_{(i)1}$ is a triideal with the even part $I_{(i)0}$ and the odd part $I_{(i)1}$ for each $i\in \Delta$. Thus, 
$Spec^{\sharp}R=\bigcup_{i\in\Delta}D^{\sharp}(I_{(i)})=
D^{\sharp}(\sum_{i\in\Delta}I_{(i)})$. Hence, $V^{\sharp}(\sum_{i\in\Delta}I_{(i)})=\emptyset$. If $1\not\in (\sum_{i\in\Delta}I_{(i)})_0=\sum_{i\in\Delta}I_{(i)0}$. Then $(\sum_{i\in\Delta}I_{(i)})_0$ is a proper ideal of the commutative ring 
$( R_0,\, +,\, \cdot)$. Hence, there exists a maximal ideal $M_0$ of the ring 
$( R_0,\, +,\, \cdot)$ such that $(\sum_{i\in\Delta}I_{(i)})_0\subseteq M_0$. Since 
$M_0\oplus R_1$ is a prime triideal of $R$ and $\sum_{i\in\Delta}I_{(i)}\subseteq M_0\oplus R_1$, we get that $M_0\oplus R_1\in V^{\sharp}(\sum_{i\in\Delta}I_{(i)})=\emptyset$, which is impossible. Therefore, $1 \in (\sum_{i\in\Delta}I_{(i)})_0=\sum_{i\in\Delta}I_{(i)0}$, which implies that $x_{(i_1)0}+x_{(i_2)0}+\cdots +x_{(i_n)0}=1$ for some positive integer $n$ and 
$x_{(i_k)0}\in I_{(i_k)0}$ with $i_k\in\Delta$ and $n\ge k\ge 1$. It follows that
$D^{\sharp}(x_{(i_k)0})\subseteq D^{\sharp}(I_{(i_k)})$ and
\begin{eqnarray*}
&&Spec^{\sharp}R=\bigcup_{i\in\Delta}D^{\sharp}(I_{(i)})\supseteq \bigcup_{k=1}^nD^{\sharp}(I_{(i_k)})\supseteq \bigcup_{k=1}^nD^{\sharp}(x_{(i_k)0})\\
&=&\bigcup_{k=1}^nD^{\sharp}(Rx_{(i_k)0})
=D^{\sharp}(\sum_{k=1}^nRx_{(i_k)0})=D^{\sharp}(R)=Spec^{\sharp}R,
\end{eqnarray*}
which implies that $Spec^{\sharp}R=\bigcup_{k=1}^nD^{\sharp}(I_{(i_k)})$.

\bigskip
(ii) Note that a closed subset of a quasicompact topological space is a quasicompact subset. Since  $Spec^{\sharp}_0R=V^{\sharp}(R_1)$ is a closed subset of the quasicompact topological space $Spec^{\sharp}R$, $Spec^{\sharp}_0R$ is a quasicompact subset.

\medskip
Note that a subset $C$ of a topological space $X$ is a quasicompact subset of $X$ if and only if every covering of $C$ by open subsets of $X$ has a finite subcovering. Hence, in order to prove that $Spec^{\sharp}_1R$ is a quasicompact subsets of $Spec^{\sharp}R$, it suffices to prove that
if $Spec^{\sharp}_1R=\bigcup_{j\in\Gamma}D^{\sharp}(J_{(j)})$ for triideals $J_{(j)}$ of $R$, then there exists a positive integer $m$ such that $Spec^{\sharp}_1R=\bigcup_{k=1}^mD^{\sharp}(J_{(j_k)})$ for some $j_1$, $\cdots$, $j_k\in\Gamma$.

\medskip
Since $Spec^{\sharp}_1R=\bigcup_{j\in\Gamma}D^{\sharp}(J_{(j)})=D^{\sharp}(\sum_{j\in\Gamma}J_{(j)})$, we have $V^{\sharp}(\sum_{j\in\Gamma}J_{(j)})=Spec^{\sharp}R\setminus D^{\sharp}(\sum_{j\in\Gamma}J_{(j)})=Spec^{\sharp}_0R$. If $(\sum_{j\in\Gamma}J_{(j)})_1\ne R_1$, then there exists a maximal ideal $N_1$ of the commutative ring 
$( R_1,\, +,\, \sharp)$ such that $(\sum_{j\in\Gamma}J_{(j)})_1\subseteq N_1$. By Proposition~\ref{pr5.1}, there exists an ideal $N_0$ of the commutative ring 
$( R_0,\, +,\, \cdot)$ such that $N_0\supseteq (\sum_{j\in\Gamma}J_{(j)})_0$ and $N_0\oplus N_1$
is a prime triideal of $R$. Thus, 
$N_0\oplus N_1\in V^{\sharp}(\sum_{j\in\Gamma}J_{(j)})=Spec^{\sharp}_0R$, which is impossible because $N_1\ne R_1$. This proves that $(\sum_{j\in\Gamma}J_{(j)})_1= R_1$. Hence, we have
$y_{(j_1)1}+y_{(j_2)1}+\cdots +y_{(j_m)1}=1^{\sharp}$ for some positive integer $m$ and $y_{(j_k)1}\in J_{(j_k)1}$ with $j_k\in\Gamma$ and $m\ge k\ge 1$. It follows that 
$D^{\sharp}(y_{(j_k)1})\subseteq D^{\sharp}(J_{(j_k)})$ and
\begin{eqnarray*}
&&Spec^{\sharp}_1R=\bigcup_{j\in\Gamma}D^{\sharp}(J_{(j)})
\supseteq \bigcup_{k=1}^mD^{\sharp}(y_{(j_k)1})\supseteq
=\supseteq \bigcup_{k=1}^mD^{\sharp}(R_1\,\sharp\,y_{(j_k)1})\supseteq\\
&=&D^{\sharp}(\sum_{k=1}^m R_1\,\sharp\,y_{(j_k)1})
=D^{\sharp}(R_1)=Spec^{\sharp}_1R,
\end{eqnarray*}
which implies that $Spec^{\sharp}_1R=\bigcup_{k=1}^mD^{\sharp}(y_{(j_k)1})$.

\bigskip
(iii) By Proposition~\ref{pr7.2},  we may assume $I=\sqrt[\sharp]{I}$ in the following proof.
First, we prove that if $V^{\sharp}(I)$ is irreducible , then $I$ is a prime triideal. 

\medskip
Suppose that $a_0b_0\in I$ for some $a_0$, $b_0\in R_0$. Let 
$$J_{(1)}=I+Ra_0=(I_0+R_0a_0)\oplus (I_1+R_1a_0)$$
and
$$K_{(1)}=I+Rb_0=(I_0+R_0b_0)\oplus (I_1+R_1b_0).$$
Then both $J_{(1)}$ and $K_{(1)}$ are triideals of $R$ and
\begin{eqnarray*}
J_{(1)}\,\stackrel{\cdot}{\sharp}\, K_{(1)}&=&(I_0+R_0a_0)(I_0+R_0b_0)+
(I_1+R_1a_0)\,\sharp\,(I_1+R_1b_0)\\
&\subseteq& I+R_0a_0R_0b_0+(R_1a_0)\,\sharp\,R_1b_0)\subseteq
I+R_0a_0b_0+R_1a_0b_0\subseteq I.
\end{eqnarray*}
Hence, we get $V^{\sharp}(J_{(1)})\cup V^{\sharp}(K_{(1)})=
V^{\sharp}(J_{(1)}\,\stackrel{\cdot}{\sharp}\, K_{(1)})\supseteq V^{\sharp}(I)$ by Proposition~\ref{pr7.1} (ii). It is clear that 
$V^{\sharp}(J_{(1)})\subseteq V^{\sharp}(I)$ and $V^{\sharp}(K_{(1)})\subseteq V^{\sharp}(I)$. Hence, we get that $V^{\sharp}(J_{(1)})\cup V^{\sharp}(K_{(1)})\subseteq  V^{\sharp}(I)$. Thus we have  $V^{\sharp}(J_{(1)})\cup V^{\sharp}(K_{(1)})= V^{\sharp}(I)$. Since $V^{\sharp}(I)$ is irreducible, $V^{\sharp}(J_{(1)})=V^{\sharp}(I)$ or $V^{\sharp}(K_{(1)})=V^{\sharp}(I)$, which imply that $I= \sqrt[\sharp]{J_{(1)}}\supseteq J_{(1)}\ni a_0$ or
$I= \sqrt[\sharp]{K_{(1)}}\supseteq K_{(1)}\ni b_0$. This proves that
\begin{equation}\label{eq52} 
a_0b_0\in I\Longrightarrow \mbox{$a_0\in I$ or $b_0\in I$ for $a_0$, $b_0\in R_0$}.
\end{equation}

\medskip
Similarly, we have
\begin{equation}\label{eq53} 
a_0b_1\in I\Longrightarrow \mbox{$a_0\in I$ or $b_1\in I$ for $a_0\in R_0$ and $b_1\in R_1$},
\end{equation}
\begin{equation}\label{eq54} 
a_1b_0\in I\Longrightarrow \mbox{$a_1\in I$ or $b_0\in I$ for $a_1\in R_1$ and $b_0\in R_0$}
\end{equation}
and
\begin{equation}\label{eq55} 
a_1\,\sharp\,b_1\in I\Longrightarrow \mbox{$a_1\in I$ or $b_1\in I$ for $a_1$, $b_1\in R_1$}.
\end{equation}

\medskip
By (\ref{eq52}), (\ref{eq53}), (\ref{eq54}) and (\ref{eq55}), $I$ is a prime triideal.

\bigskip
Next, we prove that if $I$ is a prime triideal, then $V^{\sharp}(I)$ is irreducible. Suppose that 
$V^{\sharp}(I)=V^{\sharp}(J)\cup V^{\sharp}(K)$, where $J$ and $K$ are triideals of $R$. Using Proposition~\ref{pr7.2} (ii), we can assume that $\sqrt[\sharp]{J}=J$ and $\sqrt[\sharp]{K}=K$. In this case, we have
$$
V^{\sharp}(J)\subseteq V^{\sharp}(I)\Longrightarrow I=
\sqrt[\sharp]{I}\subseteq \sqrt[\sharp]{J}=J
$$
and
$$
V^{\sharp}(K)\subseteq V^{\sharp}(I)\Longrightarrow I=
\sqrt[\sharp]{I}\subseteq \sqrt[\sharp]{K}=K.
$$
By Proposition~\ref{pr7.1} (ii), we have 
$V^{\sharp}(I)=V^{\sharp}(J)\cup V^{\sharp}(K)=V^{\sharp}(J\,\stackrel{\cdot}{\sharp}\,K)$. This fact and Proposition~\ref{pr7.2} (i) give $J\,\stackrel{\cdot}{\sharp}\,K\subseteq \sqrt[\sharp]{I}=I$. Since $I$ is a prime triideal, we get $J\subseteq I$ or $K\subseteq I$ by 
Proposition~\ref{pr5.2}. Hence, $I=J$ or $I=K$. Thus, $V^{\sharp}(I)=V^{\sharp}(J)$ or 
$V^{\sharp}(I)=V^{\sharp}(K)$. This proves that $V^{\sharp}(I)$ is irreducible.

\hfill\raisebox{1mm}{\framebox[2mm]{}}

\end{document}